\documentclass{amsproc}
\usepackage{amsthm,amsmath,amsfonts,mathrsfs,amssymb}
\title{On the centralizer of an $I$-matrix in $M_2(R/I)$,\,$I$\,a principal ideal and $R$ a UFD\quad}
\author{Magdaleen S. Marais}
\email{magdaleen@aims.ac.za}
\address{African Institute for Mathematical Sciences, 6 Melrose Rd, Muizenberg, 7945, Cape Town, South Africa }
\thanks{ This research is part of the author's research for her doctoral dissertation which was conducted at Stellenbosch University under the direction of L. van Wyk. The financial assistance of the National Research Foundation (NRF) towards this research is hereby acknowledged. Opinions expressed and conclusions arrived at are those of the author and are not necessarily to be attributed to the National Research Foundation.}
\begin{document}
\begin{abstract}
The concept of an $I$-matrix in the full $2\times 2$ matrix ring $M_2(R/I)$, where $R$ is an arbitrary UFD and $I$ is a nonzero ideal in~$R$, was introduced in~\cite{mar}. Moreover a concrete description of the centralizer of an $I$-matrix~$\widehat B$ in $M_2(R/I)$ as the sum of two subrings $\mathcal S_1$ and $\mathcal S_2$ of $M_2(R/I)$ was also given, where~$\mathcal S_1$ is the image (under the natural epimorphism from $M_2(R)$ to $M_2(R/I)$) of the centralizer in $M_2(R)$ of a pre-image of $\widehat B$, and where the entries in $\mathcal S_2$ are intersections of certain annihilators of elements arising from the entries of~$\widehat B$. In the present paper, we obtain results for the case when $I$ is a principal ideal $\langle k\rangle$, $k\in R$ a nonzero nonunit. Mainly we solve two problems. Firstly we find necessary and sufficient conditions for when $\mathcal S_1\subseteq\mathcal S_2$, for when $\mathcal S_2\subseteq \mathcal S_1$ and for when $\mathcal S_1=\mathcal S_2$. Secondly we provide a formula for the number of elements in the centralizer of $\widehat B$ for the case when $R/\langle k\rangle$ is finite.
\end{abstract}
\maketitle
\section{Introduction}
\noindent We denote the centralizer of an element $s$ in an arbitrary ring $S$ by Cen$_{S}(s)$. Knowing that $M_n(R)$, the full $n\times n$ matrix ring over a commutative ring $R$, is a prime example of a non-commutative ring, it is surprising that a concrete description of $\textnormal{Cen}_{M_n(R)}(B)$ for an arbitrary $B\in M_n(R)$ has not yet been found. If $R[x]$ is the polynomial ring in the variable $x$ over $R$, then
\begin{equation}\label{1a}
\{f(B)\ |\ f(x)\in R[x]\}\subseteq \textrm{Cen}_{M_n(R)}(B).
\end{equation}
In fact, it is known that (see \cite{cen})
\begin{equation*}
\{f(B)\ |\ f(x)\in R[x]\}=\textrm{Cen}_{M_n(R)}(\textrm{Cen}_{M_n(R)}(B)).
\end{equation*}
The most progress, finding a concrete description of $\textnormal{Cen}_{M_n(R)}(B)$, has been made for the case when the underlying ring $R$ is a field (see \cite{halm}, \cite{hung}, \cite{frob}, \cite{closed bases} and \cite{sup}).
The following well-known result in this case provides a necessary and sufficient condition for equality in (\ref{1a}).

\newtheorem{Theorem2.1}{\bf Theorem}[section]
\begin{Theorem2.1}\label{Theorem2.1}
If $B$ is an $n\times n$ matrix over a field $F$, then
\[\textnormal{Cen}_{M_n(F)}(B)=\{f(B)\ |\ f(x)\in F[x]\}\]
if and only if the minimum polynomial of $B$ coincides with the characteristic polynomial of
$B$.
\end{Theorem2.1}
 The concept of an $I$-matrix in the full $2\times 2$ matrix ring $M_2(R/I)$, where $R$ is a UFD and $I$ an ideal in $R$ was introduced in \cite{mar}. In this paper, unless stated otherwise, we assume that $R$ is a UFD, $I$ is a nonzero ideal in $R$ and~$k:=\textnormal{gcd}(I)\neq 0$.
Let $\theta_I:R\to R/I$ and $\Theta_I:M_2(R)\to
M_2(R/I)$ be the
natural epimorphism and induced epimorphism respectively. We denote the image $\theta_I(b)$ of~$b\in R$ by $\hat b_I$ and the image $\Theta_I(B)$ of $B\in M_2(R)$ by $\widehat B_I$. However, if there is no ambiguity, then we simply write $\theta$, $\Theta$,~$\hat b$ and $\widehat B$ respectively.

\newtheorem{Def10a}{\bf Definition}[section]
\begin{Def10a}\label{Def10a}
We call a matrix $\left[\begin{array}{cc}\hat e_I&\hat f_I\\\hat g_I&\hat h_I\end{array}\right]\in M_2(R/I)$ an $I$-matrix \linebreak if $\langle \hat e_I-\hat h_I,\hat f_I\rangle=\langle\hat t_I\rangle$ or $\langle \hat e_I-\hat h_I, \hat g_I\rangle=\langle\hat t_I\rangle$ or $\langle \hat f_I, \hat g_I\rangle=\langle\hat t_I\rangle$, where $t|k$.
\end{Def10a}

If $R$ is a PID, then every matrix in $M_2(R/I)$ is an $I$-matrix.
A concrete description of the centralizer of an $I$-matrix, as the sum of two subrings of $M_2(R/I)$, was given by the following result in \cite{mar}:
\newtheorem{Theorem21}[Def10a]{\bf Theorem}
\begin{Theorem21}\label{Theorem21}
Let $R$ be a UFD, $I$ a nonzero ideal in $R$, and let $\widehat B_I=\left[\begin{array}{cc}\hat e_I&\hat f_I\\\hat g_I&\hat h_I\end{array}\right]\in M_2(R/I)$ be an $I$-matrix, then $\textnormal{Cen}(\widehat B)
=\mathcal S_1+\mathcal S_2$,
where\newline \[\mathcal S_1=\Theta(\textnormal{Cen}(B))\quad\textnormal{ and }\quad\mathcal S_2=\left[\begin{array}{cc}\textnormal{ann}(\hat f)\cap \textnormal{ann}(\hat g)&\textnormal{ann}(\hat g)\cap
\textnormal{ann}(\hat e-\hat h)\\\textnormal{ann}(\hat f)\cap
\textnormal{ann}(\hat e-\hat h)&\textnormal{ann}(\hat f)\cap \textnormal{ann}(\hat g)\end{array}\right].\]
\end{Theorem21}
Unfortunately the concrete description in Theorem \ref{Theorem21} could not be generalized to $n\times n$-matrices, for $n\ge 3$, in the sense of Proposition \ref{Lemma5}. In \cite{mar}, for $R$ a UFD, a matrix was given for every factor ring $R/I$ with zero divisors and every $n\ge 3$ for which equality in (\ref{verg}) does not hold.
\newtheorem{Lemma5}[Theorem2.1]{\bf Proposition}
\begin{Lemma5}\label{Lemma5}
Let $R$ be a commutative ring and let $B=[b_{ij}]\in M_n(R)$. Then
\begin{equation}\label{verg}\Theta(\textnormal{Cen}(B))+[\mathcal{A}_{ij}]\subseteq
\textnormal{Cen}(\widehat B),\end{equation} where
\[\mathcal{A}_{ij}=\left(\displaystyle\bigcap_{k,\ k\neq j} \textnormal{ann}(\hat b_{jk})\right)\bigcap\left(\displaystyle\bigcap_{k,\ k\neq i} \textnormal{ann}(\hat b_{ki})\right)\bigcap\ \textnormal{ann}(\hat b_{ii}-\hat b_{jj}).\]
\end{Lemma5}

Regarding Theorem \ref{Theorem21}, an example was also provided in \cite{mar}, where $\mathcal S_1\not\subseteq\mathcal S_2$ and $\mathcal S_2\not\subseteq\mathcal S_1$, from which the following questions arise: When is $\mathcal S_1\not\subseteq\mathcal S_2$, when is $\mathcal S_2\not\subseteq\mathcal S_1$ and when is $\mathcal S_1=\mathcal S_2$? In Section \ref{containment} this questions will be answered for the case when $I\subset R$ is a principal ideal $\langle k\rangle$ generated by a nonzero nonunit $k\in R$.

The problem of enumerating the number of matrices with given characteristics over a finite ring has been treated extensively in the literature. Formulas have been found, for example, for the number of matrices with a given characteristic polynomial \cite{char}; the number of matrices over a finite field that are cyclic \cite{ber} or symmetric \cite{car}; and the number of matrices over the ring of integers $\mathbb Z$ modulo $m$, $\mathbb Z_m$, that are nilpotent \cite{bol1}. By using the results in \cite{bol}, some of the above mentioned results, where the matrices over a finite field that satisfy some property are enumerated by rank, can be extended to matrices over certain finite rings that satisfy the property under consideration.

Naturally the question whether it is possible to enumerate the number of matrices in $\textnormal{Cen}_{M_n(R)}(B)$, denoted by $|\textnormal{Cen}_{M_n(R)}(B)|$, when~$R$ is a finite commutative ring and $B\in M_n(R)$, arises. Using the fact that the dimension of $\textnormal{Cen}_{M_n(F)}(B)$ is known by the following result, due to Frobenius, the answer is straightforward in the case when $R$ is a finite field $F$.

\newtheorem{Frobenius}[Def10a]{\bf Theorem}
\begin{Frobenius}\label{Frobenius}
Let $B\in M_n(F)$, and suppose that $f_1,\ldots,f_l\in F[x]$ are the invariant factors of $B$, where $f_i$ divides $f_{i-1}$, for~$i=2,\ldots,l$. Then the dimension of $Cen_{M_n(F)}(B)$ is given by
\[\sum_{i=1}^{l}(2i-1)(\deg f_i).\]
\end{Frobenius}
For example, if $n=2$, then the number of elements in $\textnormal{Cen}_{M_n(F)}(B)$ is~$|F|^2$, if~$B$ is a nonscalar matrix, and it is $|F|^4$ if~$B$ is a scalar matrix. Unfortunately the answer is not that easy in the case when $R$ is a finite ring that is not a field.

In Section~\ref{finite} we define an equivalence relation on $M_2(R/\langle k\rangle)$ and we use this relation, together with Theorem \ref{Theorem21}, and the results in Section \ref{containment}, to obtain a formula for the number of matrices in $\textnormal{Cen}_{M_2(R/\langle k\rangle)}(\widehat B)$ when $R$ is a UFD and $R/\langle k\rangle$ is finite, $k$ is a nonzero nonunit element in $R$ and $\widehat B\in M_2(R/\langle k\rangle)$. 

\section{Containment considerations regarding the centralizer of a $\langle k\rangle$-matrix}\label{containment}
In this section we answer the following questions: Regarding Theorem \ref{Theorem21}, when is $\mathcal S_1\not\subseteq\mathcal S_2$, when is $\mathcal S_2\not\subseteq\mathcal S_1$ and when is $\mathcal S_1=\mathcal S_2$?

We need the following preliminary definitions and results in Theorem \ref{Corollary12c}, the main result of this section.

Since the minimum polynomial and characteristic polynomial of any $2\times 2$ non-scalar matrix over a field coincide, Theorem~\ref{Theorem2.1} can be written in the following form for the $2\times 2$ case.
\newtheorem{Corollary2.5}[Theorem2.1]{\bf Corollary}
\begin{Corollary2.5}\label{Corollary2.5}
Let $B=\left[\begin{array}{cc}e&f\\g&h\end{array}\right]\in M_2(F)$, $F$ a field. Then \[
{\textnormal{Cen}_{M_2(F)}(B)
=\left\{\begin{array}{l}(i)\:M_2(F),\ \textit{if $e=h$, $f=0$ and $g=0$ (i.e.~$B$ is a
scalar matrix)}\\\\
(ii)\left.\left\{\left[\begin{array}{cc}a&0\\0&b\end{array}\right]\right|a,b\in
F\right\},\ \textit{if
$e\neq h$, $f=0$ and $g=0$}\\\\
 (iii)\left.\left\{\left[\begin{array}{cc}a&0\\b&a-g^{-1}b(e-h)\end{array}\right]\right|a,b\in
 F\right\},\ \textit{if $f=0$ and $g\neq 0$}\\\\
 (iv){\left.\left\{\left[\begin{array}{cc}a&b\\f^{-1}gb&a-f^{-1}b(e-h)\end{array}\right]\right|a,b\in
 F\right\},\ \textit{if $f\neq 0$}}.\end{array}\right.}\]
\end{Corollary2.5}

The following result, giving a concrete description of the centralizer of a matrix in $M_2(R)$, was proved in \cite{mar}:

\newtheorem{Corollary2.5b}[Theorem2.1]{\bf Lemma}
\begin{Corollary2.5b}\label{Corollary2.5b}
Let $B=\left[\begin{array}{cc}e&f\\g&h\end{array}\right]\in M_2(R)$, $R$ a UFD. Then
$\textnormal{Cen}_{M_2(R)}(B)$
\[=\left\{\begin{array}{l}(i)\:M_2(R),\ \textit{if $e=h$, $f=0$ and $g=0$ (i.e.~$B$ is a
scalar matrix)}\\\\
(ii)\left.\left\{m^{-1}w\left[\begin{array}{cc}e-h&f\\g&0\end{array}\right]+vE\right|v,w\in
 R\right\}, \begin{array}{l}\textit{if at least one}\\\textit{of $e-h,f,g$ is nonzero,}\end{array}\end{array}\right.\]
 where $m^{-1}$ is the inverse of $m:=\textnormal{gcd}(e-h,f,g)$ in the quotient field of~$R$.
\end{Corollary2.5b}

The following four results can be easily proved.

\newtheorem{newlemma}[Theorem2.1]{\bf Lemma}
\begin{newlemma}\label{newlemma}
Let $S$ be a subring of a ring $T$ and let $s\in S$.
Then
\[\textnormal{Cen}_S(s)=S\cap \textnormal{Cen}_T(s).\]
\end{newlemma}
\newtheorem{Theorem1.1}[Theorem2.1]{\bf Lemma}
\begin{Theorem1.1}\label{Theorem1.1:box}
Let $B\in M_n(R)$, where $R$ is a commutative ring. Then
\[\textnormal{Cen}_{M_2(R)}(B^T)=(\textnormal{Cen}_{M_2(R)}(B))^T.\]
\end{Theorem1.1}

\newtheorem{Lemma11}[Theorem2.1]{\bf Lemma}
\begin{Lemma11}\label{Lemma11}
Let $R$ be a UFD. Suppose $b,k\in R$, $k$ a nonzero nonunit, and $\delta=\textnormal{gcd}(b,k)$. Then
\[\langle t\rangle=\theta^{-1}(\textnormal{ann}(\hat
b_{\langle k\rangle})),\] where $t=\delta^{-1}k\in R$, with $\delta^{-1}$ the inverse of $\delta$ in the quotient field
of $R$.
\end{Lemma11}

\newtheorem{Lemma21}[Theorem2.1]{\bf Lemma}
\begin{Lemma21}\label{Lemma21}
Let $R$ be a UFD and let $k,x,y\in R$, then
\[\textnormal{ann}(\hat d)=\textnormal{ann}(\hat x)\cap \textnormal{ann}(\hat y)\]
in $R/\langle k\rangle$, with \textnormal{gcd}$(x,y)=d$.
\end{Lemma21}

We are now in a position to prove Theorem \ref{Lemma21}.

\newtheorem{Corollary12c}[Theorem2.1]{\bf Theorem}
\begin{Corollary12c}\label{Corollary12c}
Let $R$ be a UFD, $k=p_1^{n_1}p_2^{n_2}\cdots p_m^{n_m}$ and let\newline $B=\left[\begin{array}{cc}e&f\\g&h\end{array}\right]\in M_2(R)$ be such that $\widehat B$ is a $\langle k\rangle$-matrix. Then
\begin{itemize}
\item[(a)]\begin{equation}\label{centr}
\textnormal{Cen}_{M_2(R/\langle k\rangle)}(\widehat B)=\Theta_{\langle k\rangle}(\textnormal{Cen}_{M_2(R)}(B))
\end{equation}
if and only if $B$ is a scalar matrix or satisfies the following conditions for every $i$, $i=1,2,\ldots,m$:
\begin{itemize}
\item[(i)] $p_i$ is not a divisor of at least one of the elements $e-h$, $f$ and $g$; pick such an element $a$, and call the remaining two elements $b$ and $c$, say.
\item[(ii)] $\textnormal{gcd}(b,c,k)=1$ or $\hat a_{\langle\textnormal{gcd}(b,c,k)\rangle}$ is invertible in $R/\langle \textnormal{gcd}(b,c,k)\rangle$;
\end{itemize}
\item[(b)]\begin{equation}\label{bevat1}
\textnormal{Cen}(\widehat B)=\left[\begin{array}{cc}\textnormal{ann}(\hat f)\cap \textnormal{ann}(\hat g)&\textnormal{ann}(\hat g)\cap \textnormal{ann}(\hat e-\hat h)\\\textnormal{ann}(\hat f)\cap \textnormal{ann}(\hat e-\hat h)&\textnormal{ann}(\hat f)\cap \textnormal{ann}(\hat g)\end{array}\right]\end{equation}
if and only if $\hat f=\hat 0$ and $\hat g=\hat 0$;

\item[(c)]\begin{equation}\label{bevat2}\Theta(\textnormal{Cen}(B))=\left[\begin{array}{cc}\textnormal{ann}(\hat f)\cap \textnormal{ann}(\hat g)&\textnormal{ann}(\hat g)\cap \textnormal{ann}(\hat e-\hat h)\\\textnormal{ann}(\hat f)\cap \textnormal{ann}(\hat e-\hat h)&\textnormal{ann}(\hat f)\cap \textnormal{ann}(\hat g)\end{array}\right]\end{equation} if and only if $\hat f=\hat 0,$ $\hat g=\hat 0$ and ($\hat e-\hat h$ is invertible or $\hat e-\hat h=\hat 0$).

\end{itemize}
\end{Corollary12c}

\begin{proof}
\noindent(a) Since (\ref{centr}) follows trivially if $B$ is a scalar matrix, we assume that $B$ is a non-scalar matrix. Suppose that conditions (i) and (ii) are satisfied. If 
\begin{equation}\label{centr1}
\textnormal{ann}_{M_2(R/\langle k\rangle)}(\theta_{\langle k\rangle}(\textnormal{gcd}(f,g)))=\hat 0_{\langle k\rangle},\qquad \textnormal{ann}_{M_2(R/\langle k\rangle)}(\theta_{\langle k\rangle}(\textnormal{gcd}(f,e-h)))=\hat 0_{\langle k\rangle}\end{equation}\begin{equation}\label{centr2}\textrm{and}\qquad \textnormal{ann}_{M_2(R/\langle k\rangle)}(\theta_{\langle k\rangle}(\textnormal{gcd}(g,e-h)))=\hat 0_{\langle k\rangle},
\end{equation}
then the result follows from Theorem \ref{Theorem21} and Lemma \ref{Lemma21}. Thus suppose that at least one of the annihilators in (\ref{centr1}) and (\ref{centr2}) is nonzero. We now show that
\begin{equation*}
\left[\begin{array}{cc}\hat 0_{\langle k\rangle}&\textnormal{ann}(\theta_{\langle k\rangle}(\textnormal{gcd}(g,e-h)))\\\hat 0_{\langle k\rangle}&\hat 0_{\langle k\rangle}\end{array}\right],\left[\begin{array}{cc}\hat 0_{\langle k\rangle}&\hat 0_{\langle k\rangle}\\\textnormal{ann}(\theta_{\langle k\rangle}(\textnormal{gcd}(f,e-h)))&\hat 0_{\langle k\rangle}\end{array}\right],
\end{equation*}
\begin{equation}\label{centr3}
\left[\begin{array}{cc}\hat 0_{\langle k\rangle}&\hat 0_{\langle k\rangle}\\\hat 0_{\langle k\rangle}&\textnormal{ann}(\theta_{\langle k\rangle}(\textnormal{gcd}(f,g)))\end{array}\right]\in\Theta_{\langle k\rangle}(\textnormal{Cen}_{M_2(R)}(B)).
\end{equation}
Since then, because $\Theta_{\langle k\rangle}(\textnormal{Cen}_{M_2(R)}(B))$ is a ring, (\ref{centr}) follows from Theorem \ref{Theorem21} and Lemma~\ref{Lemma21}.

If $\textnormal{ann}_{M_2(R/\langle k\rangle)}(\theta_{\langle k\rangle}(\textnormal{gcd}(g,e-h)))\neq\hat 0_{\langle k\rangle}$, then, by Lemma \ref{Lemma11}, \newline$1\neq\textnormal{gcd}(g,e-h,k):=\delta$ and $\textnormal{ann}_{M_2(R/\langle k\rangle)}(\theta_{\langle k\rangle}(\textnormal{gcd}(e-h,g)))=\langle (\widehat{k\delta^{-1}})_{\langle k\rangle}\rangle$. To accomplish our objective, we show that for each $\hat d_{\langle k\rangle}\in \textnormal{ann}_{M_2(R/\langle k\rangle)}(\theta_{\langle k\rangle}(\textnormal{gcd}(g,e-h)))$ there is a $\hat d'_{\langle k\rangle}\in \textnormal{ann}_{M_2(R/\langle k\rangle)}(\theta_{\langle k\rangle}(\textnormal{gcd}(g,e-h)))$ such that $\hat f_{\langle k\rangle}\hat d'_{\langle k\rangle}=\hat d_{\langle k\rangle}$, since then
\[\Theta_{\langle k\rangle}\left(\left[\begin{array}{cc}0&fd'\\gd'&(e-h)d'\end{array}\right]\right)=\left[\begin{array}{cc}\hat 0_{\langle k\rangle}&\hat d_{\langle k\rangle}\\\hat 0_{\langle k\rangle}&\hat 0_{\langle k\rangle}\end{array}\right],\]
so that we therefore can conclude from Lemma \ref{Corollary2.5b}(ii) that
\[\left[\begin{array}{cc}\hat 0_{\langle k\rangle}&\textnormal{ann}_{M_2(R/\langle k\rangle)}(\theta_{\langle k\rangle}(\textnormal{gcd}(g,e-h)))\\\hat 0_{\langle k\rangle}&\hat 0_{\langle k\rangle}\end{array}\right]\in\Theta_{\langle k\rangle}(\textnormal{Cen}_{M_2(R)}(B)).\]
Thus, let $\hat d_{\langle k\rangle}$ be an arbitrary element in $\textnormal{ann}_{M_2(R/\langle k\rangle)}(\theta_{\langle k\rangle}(\textnormal{gcd}(g,e-h)))$, i.e.~suppose~$\hat d_{\langle k\rangle}:=\hat s_{\langle k\rangle}(\widehat{k\delta^{-1}})_{\langle k\rangle}$ for some $\hat s_{\langle k\rangle}\in R/\langle k\rangle$. Since $\hat f_{\langle\delta\rangle}$ is invertible in $R/\langle \delta\rangle$, by (ii), there is a $\hat t_{\langle\delta\rangle}\in R/\langle\delta\rangle$ such that $\hat t_{\langle\delta\rangle}\hat f_{\langle\delta\rangle}=\hat 1_{\langle\delta\rangle}$ which implies that $tf=1+v\delta$ for some $v\in R$. Hence $ftd=(1+v\delta)(sk\delta^{-1}+wk)=sk\delta^{-1}+(w+vs+v\delta w)k$. In other words, if we set $\hat d'_{\langle k\rangle}:=(\widehat{td})_{\langle k\rangle}$ then $\hat f_{\langle k\rangle}\hat d'_{\langle k\rangle}=\hat f_{\langle k\rangle}(\widehat{td})_{\langle k\rangle}=(\widehat{sk\delta^{-1}})_{\langle k\rangle}=\hat d_{\langle k\rangle}$.

It can similarly be shown that each of the other two sets in (\ref{centr3}) is contained in $\Theta_{\langle k\rangle}(\textnormal{Cen}_{M_2(R)}(B))$.\vskip 0.5 cm

Conversely, suppose $B$ does not satisfy both of the conditions (i) and (ii). We distinguish between the following cases:
\begin{itemize}
\item[(a$'$)] $B$ does not satisfy (i), i.e.~$\textnormal{gcd}(e-h,f,g,k)\neq 1$;
\item[(b$'$)] $B$ satisfies (i), but not (ii).
\end{itemize}\vskip 0.5cm

\noindent(a$'$) Suppose there is a prime $p_i$ in the prime factorization of $k$ such that $p_i|e-h,f,g$. We distinguish between the following two cases:
\begin{itemize}
\item[(i$'$)]$f=0$ or $g=0$;
\item[(ii$'$)]$f,g\neq0$.
\end{itemize}\vskip 0.3cm

\noindent(i$'$) Since $p_i|e-h,f,g$, direct verification shows that \begin{eqnarray*}\widehat A_{\langle k\rangle}&:=&\left[\begin{array}{cc}\hat 0&\theta(p_1^{n_1}\cdots p_{i-1}^{n_{i-1}}p_i^{n_i-1}p_{i+1}^{n_{i+1}}\cdots p_m^{n_m})\\\theta(p_1^{n_1}\cdots p_{i-1}^{n_{i-1}}p_i^{n_i-1}p_{i+1}^{n_{i+1}}\cdots p_m^{n_m})&\hat{0}\end{array}\right]\\&\in&\textnormal{Cen}_{M_2(R/\langle k\rangle)}(\widehat B_{\langle k\rangle}).\end{eqnarray*}
Because $\theta_{\langle k\rangle}(p_1^{n_1}\cdots p_{i-1}^{n_{i-1}}p_i^{n_i-1}p_{i+1}^{n_{i+1}}\cdots p_m^{n_m})\neq\hat 0_{\langle k\rangle}$, it follows that the entries in position~$(1,2)$ and position~$(2,1)$ of $\widehat A_{\langle k\rangle}$ only have nonzero pre-images in $R$. Since~$B$ is a non-scalar matrix, it follows from Lemma \ref{Corollary2.5b}(ii) that
every matrix in $\textnormal{Cen}_{M_2(R)}(B)$ has $0$ in position $(1,2)$ if $f=0$ and $0$ in position $(2,1)$ if $g=0$. Therefore $\widehat A_{\langle k\rangle}\not\in \Theta_{\langle k\rangle}(\textnormal{Cen}_{M_2(R)}(B))$ if $f=0$ or $g=0$.\vskip 0.3cm

\noindent(ii$'$) Since $f,g\neq0$ and $p_i|f,g$ it follows that 
\begin{equation}\label{centr4}
f=cp^r_i\qquad\textrm{ and }\qquad g=dp^s_i
\end{equation}
for some $s,r\ge1$ and $c,d\in R$ such that $p_i\nmid c,d$. Now, $r\le s$ or $s\le r$. Let us first assume that $r\le s$. Because $p_i|e-h,f,g$ direct verification shows that
\[\widehat A_{\langle k\rangle}:=\left[\begin{array}{cc}\hat 0_{\langle k\rangle}&\hat 0_{\langle k\rangle}\\\theta_{\langle k\rangle}(p_1^{n_1}\cdots p_{i-1}^{n_{i-1}}p_i^{n_i-1}p_{i+1}^{n_{i+1}}\cdots p_m^{n_m})&\hat 0_{\langle k\rangle}\end{array}\right]\in \textnormal{Cen}_{M_2(R/\langle k\rangle)}(\widehat B_{\langle k\rangle}).\] 

\noindent We now show that $\widehat A_{\langle k\rangle}\not\in\Theta_{\langle k\rangle}(\textnormal{Cen}_{M_2(R)}(B))$. 
Firstly note that the set of all the pre-images of~$\widehat A_{\langle k\rangle}$ is
\[\left[\begin{array}{cc}\ker\theta_{\langle k\rangle}&\ker\theta_{\langle k\rangle}\\p_1^{n_1}\cdots p_{i-1}^{n_{i-1}}p_i^{n_i-1}p_{i+1}^{n_{i+1}}\cdots p_m^{n_m}+\ker\theta_{\langle k\rangle}&\ker\theta_{\langle k\rangle}\end{array}\right].\]Thus, if $\widehat A_{\langle k\rangle}\in\Theta_{\langle k\rangle}(\textnormal{Cen}_{M_2(R)}(B))$, then, according to Corollary~\ref{Corollary2.5}(iv) and Lem-ma~\ref{newlemma} there is a pre-image $\left[\begin{array}{cc}\kappa_1&\kappa_2\\p_1^{n_1}\cdots p_{i-1}^{n_{i-1}}p_i^{n_i-1}p_{i+1}^{n_{i+1}}\cdots p_m^{n_m}+\kappa_3&\kappa_4\end{array}\right]\in M_2(R)$ of $\widehat A_{\langle k\rangle}$, where $\kappa_1,\kappa_2,\kappa_3,\kappa_4\in\ker\theta_{\langle k\rangle}$, such that
\[\left[\begin{array}{cc}\kappa_1&\kappa_2\\p_1^{n_1}\cdots p_{i-1}^{n_{i-1}}p_i^{n_i-1}p_{i+1}^{n_{i+1}}\cdots p_m^{n_m}+\kappa_3&\kappa_4\end{array}\right]=\left[\begin{array}{cc}a&b\\gf^{-1}b&a-(e-h)f^{-1}b\end{array}\right]\]in $M_2(R)$ for some $a,b\in R$. In other words, there are $a,b\in R$ such that $\kappa_1=a$, $\kappa_2=b$ and $p_1^{n_1}\cdots p_{i-1}^{n_{i-1}}p_i^{n_i-1}p_{i+1}^{n_{i+1}}\cdots p_m^{n_m}+\kappa_3=gf^{-1}b$. But then, considering~(\ref{centr4}) and keeping in mind that $r\le s$, $gf^{-1}b\in R$ and $p_i\nmid c,d$, we have that $gf^{-1}b=\nolinebreak dp_i^s(cp_i^r)^{-1}\kappa_2\in\langle p_i^{n_i}\rangle$, where $\langle p_i^{n_i}\rangle$ is the ideal generated by $p_i^{n_i}$ in $R$. Because $p_i^{n_i}\nmid p_1^{n_1}\cdots p_{i-1}^{n_{i-1}}p_i^{n_i-1}p_{i+1}^{n_{i+1}}\cdots p_m^{n_m}+\kappa_3$, it follows that $p_1^{n_1}\cdots p_{i-1}^{n_{i-1}}p_i^{n_i-1}p_{i+1}^{n_{i+1}}$ $\cdots p_m^{n_m}+\kappa_3\not\in\langle p_i^{n_i}\rangle$, which implies that \[p_1^{n_1}\cdots p_{i-1}^{n_{i-1}}p_i^{n_i-1}p_{i+1}^{n_{i+1}}\cdots p_m^{n_m}+\kappa_3\neq gf^{-1}b.\] Thus we have a contradiction. Therefore $\widehat A_{\langle k\rangle}\not\in\Theta_{\langle k\rangle}(\textnormal{Cen}_{M_2(R)}(B))$. 

If $s\le r$ one can similarly show that \[\widehat A_{\langle k\rangle}:=\left[\begin{array}{cc}\hat 0_{\langle k\rangle}&\theta_{\langle k\rangle}(p_1^{n_1}\cdots p_{i-1}^{n_{i-1}}p_i^{n_i-1}p_{i+1}^{n_{i+1}}\cdots p_m^{n_m})\\\hat 0_{\langle k\rangle}&\hat 0_{\langle k\rangle}\end{array}\right]\in \textnormal{Cen}_{M_2(R/\langle k\rangle)}(\widehat B_{\langle k\rangle}),\] and that $\widehat A_{\langle k\rangle}\not\in\Theta_{\langle k\rangle}(\textnormal{Cen}_{M_2(R)}(B))$, by using Lemma \ref{Theorem1.1:box} and Corollary \ref{Corollary2.5}(iv) instead of Corollary \ref{Corollary2.5}(iv).\vskip 0.5cm

\noindent (b$'$) Suppose $B$ satisfies (i), but not (ii). Then, for each $i\in\{1,\ldots,m\}$, at least one of the following cases is true:
\begin{itemize}
\item[(i$'$)] $\textnormal{gcd}(e-h,f,g,k)=1$ , $1\neq\textnormal{gcd}(e-h,g,k):=\delta$ and $\hat f_{\langle\delta\rangle}$ is not invertible in~$R/\langle \delta\rangle$;
\item[(ii$'$)] $\textnormal{gcd}(e-h,f,g,k)=1$ , $1\neq\textnormal{gcd}(e-h,f,k):=\delta$ and $\hat g_{\langle\delta\rangle}$ is not invertible in~$R/\langle \delta\rangle$;
\item[(iii$'$)] $\textnormal{gcd}(e-h,f,g,k)=1$ , $1\neq\textnormal{gcd}(f,g,k):=\delta$ and $\hat e_{\langle\delta\rangle}-\hat h_{\langle\delta\rangle}$ is not invertible in $R/\langle \delta\rangle$;
\end{itemize}
We now show that (\ref{centr}) does not follow in each of the above cases.\vskip 0.3cm

\noindent (i$'$) In this case Lemma \ref{Lemma11} implies that \[\textnormal{ann}_{M_2(R/\langle k\rangle)}(\theta_{\langle k\rangle}(\textnormal{gcd}(g,e-h)))=\langle(\widehat{k\delta^{-1}})_{\langle k\rangle}\rangle.\] Note that since $\delta$ is not a unit, $\langle k\delta^{-1}\rangle\neq \langle k\rangle$. By Theorem \ref{Theorem21} it follows that
\[\widehat A_{\langle k\rangle}:=\left[\begin{array}{cc}\hat 0_{\langle k\rangle}&(\widehat{k\delta^{-1}})_{\langle k\rangle}\\\hat 0_{\langle k\rangle}&\hat 0_{\langle k\rangle}\end{array}\right]\in\textnormal{Cen}_{M_2(R/\langle k\rangle)}(\widehat B).\]
If we can show that $\widehat A_{\langle k\rangle}\not\in\Theta_{\langle k\rangle}(\textnormal{Cen}_{M_2(R)}(B))$, then we are finished. Now,
\[\left[\begin{array}{cc}\ker\theta_{\langle k\rangle}&k^{-1}\delta+\ker\theta_{\langle k\rangle}\\\ker\theta_{\langle k\rangle}&\ker\theta_{\langle k\rangle}\end{array}\right]\] is the set of all the pre-images of $\widehat A_{\langle k\rangle}$ in $R$. Therefore, taking into account that $\textnormal{gcd}(e-h,f,g,k)=1$, if $\widehat A_{\langle k\rangle}\in\Theta_{\langle k\rangle}(\textnormal{Cen}_{M_2(R)}(B))$ it follows from Corollary \ref{Corollary2.5b}(ii) that there is a pre-image $\left[\begin{array}{cc}\kappa_1&k\delta^{-1}+\kappa_2\\\kappa_3&\kappa_4\end{array}\right]\in M_2(R)$ of $\widehat A_{\langle k\rangle}$, where $\kappa_1, \kappa_2, \kappa_3, \kappa_4\in\ker\theta_{\langle k\rangle}$, such that 
\[\left[\begin{array}{cc}\kappa_1&k\delta^{-1}+\kappa_2\\\kappa_3&\kappa_4\end{array}\right]=\left[\begin{array}{cc}a&fb\\gb&a-(e-h)b\end{array}\right]\] for some $a,b\in R$. Hence, $gb=\kappa_3$ and $(e-h)b=\kappa_1-\kappa_4$, which implies that $b=sk\delta^{-1}$ for some $s\in R$. But then, since $fb=k\delta^{-1}+\kappa_2$, we have that
\begin{equation*}fb=fsk\delta^{-1}=k\delta^{-1}+\kappa_2\Leftrightarrow fs=1+t\delta\ \textrm{for some }t\in R\Leftrightarrow \hat f_{\langle\delta\rangle}\hat s_{\langle\delta\rangle}=\hat 1_{\langle\delta\rangle}.\end{equation*}Since $\hat f_{\langle\delta\rangle}$ is not invertible in $R/\langle\delta\rangle$, according to assumption, we have a contradiction. Therefore $\widehat A_{\langle k\rangle}\not\in\Theta_{\langle k\rangle}(\textnormal{Cen}_{M_2(R)}(B))$ and so we conclude that \[\textnormal{Cen}_{M_2(R/\langle k\rangle)}(\widehat B_{\langle k\rangle})\not\subseteq\Theta_{\langle k\rangle}(\textnormal{Cen}_{M_2(R)}(B)).\]
\vskip 0.3cm
\noindent(ii$'$ and iii$'$) It follows similarly to case (i$'$) that $\textnormal{Cen}_{M_2(R/\langle k\rangle)}(\widehat B)\not\subseteq\Theta(\textnormal{Cen}_{M_2(R)}(B))$.\vskip 0.5cm

\noindent (b) Suppose $\hat f,\hat g=\hat 0$, then $f,g\in\langle k\rangle$, and so by Corollary \ref{Corollary2.5b}(ii)
\begin{eqnarray*}
\Theta(\textnormal{Cen}(B))&=&\Theta\left(\left.\left\{\left[\begin{array}{cc}a&fb\\gb&a-(e-h)b\end{array}\right]\right|a,b\in R\right\}\right)\\&=&\Theta\left(\left.\left\{\left[\begin{array}{cc}a&0\\0&a-(e-h)b\end{array}\right]\right|a,b\in R\right\}\right)\\
&\subseteq&\left[\begin{array}{cc}R/\langle k\rangle&\textnormal{ann}(\hat e-\hat h)\\\textnormal{ann}(\hat e-\hat h)&R/\langle k\rangle\end{array}\right]\\&=&\left[\begin{array}{cc}\textnormal{ann}(\hat f)\cap\textnormal{ann}(\hat g)&\textnormal{ann}(\hat g)\cap\textnormal{ann}(\hat e-\hat h)\\\textnormal{ann}(\hat f)\cap\textnormal{ann}(\hat e-\hat h)&\textnormal{ann}(\hat f)\cap\textnormal{ann}(\hat g)\end{array}\right].\end{eqnarray*}
Conversely, suppose $\Theta(\textnormal{Cen}(B))\subseteq \left[\begin{array}{cc}\textnormal{ann}(\hat f)\cap\textnormal{ann}(\hat g)&\textnormal{ann}(\hat g)\cap\textnormal{ann}(\hat e-\hat h)\\\textnormal{ann}(\hat f)\cap\textnormal{ann}(\hat e-\hat h)&\textnormal{ann}(\hat f)\cap\textnormal{ann}(\hat g)\end{array}\right]$. Since $\left[\begin{array}{cc}\hat a&\hat 0\\\hat 0&\hat a\end{array}\right]\in\Theta(\textnormal{Cen}(B))$ for every $\hat a\in R/\langle k\rangle$ it follows that \newline$\textnormal{ann}(\hat f)\cap\textnormal{ann}(\hat g)=R/\langle k\rangle$ which implies that $\textnormal{ann}(\hat f)=R/\langle k\rangle$ and $\textnormal{ann}(\hat g)=R/\langle k\rangle$ and so $\hat f,\hat g=\hat 0$.\vskip 0.15cm

\noindent(c) Using (b) and (a), it follows that \begin{eqnarray*}
&&\Theta(\textnormal{Cen}(B))=\left[\begin{array}{cc}\textnormal{ann}(\hat f)\cap\textnormal{ann}(\hat g)&\textnormal{ann}(\hat g)\cap\textnormal{ann}(\hat e-\hat h)\\\textnormal{ann}(\hat f)\cap\textnormal{ann}(\hat e-\hat h)&\textnormal{ann}(\hat f)\cap\textnormal{ann}(\hat g)\end{array}\right]\nonumber\\
&\Leftrightarrow& \Theta(\textnormal{Cen}(B))\subseteq\left[\begin{array}{cc}\textnormal{ann}(\hat f)\cap\textnormal{ann}(\hat g)&\textnormal{ann}(\hat g)\cap\textnormal{ann}(\hat e-\hat h)\\\textnormal{ann}(\hat f)\cap\textnormal{ann}(\hat e-\hat h)&\textnormal{ann}(\hat f)\cap\textnormal{ann}(\hat g)\end{array}\right]\textrm{ and }\nonumber\\&&\left[\begin{array}{cc}\textnormal{ann}(\hat f)\cap\textnormal{ann}(\hat g)&\textnormal{ann}(\hat g)\cap\textnormal{ann}(\hat e-\hat h)\\\textnormal{ann}(\hat f)\cap\textnormal{ann}(\hat e-\hat h)&\textnormal{ann}(\hat f)\cap\textnormal{ann}(\hat g)\end{array}\right]\subseteq\Theta(\textnormal{Cen}(B))\nonumber\\
&\Leftrightarrow&\hat f,\hat g=\hat 0 \textrm{ and }\left[\begin{array}{cc}\textnormal{ann}(\hat f)\cap\textnormal{ann}(\hat g)&\textnormal{ann}(\hat g)\cap\textnormal{ann}(\hat e-\hat h)\\\textnormal{ann}(\hat f)\cap\textnormal{ann}(\hat e-\hat h)&\textnormal{ann}(\hat f)\cap\textnormal{ann}(\hat g)\end{array}\right]\subseteq\Theta(\textnormal{Cen}(B))\label{c1}\\
&\Leftrightarrow&\hat f,\hat g=\hat 0\textrm{ and }(\hat e-\hat h \textrm{ is invertible in }R/\langle k\rangle\textrm{ or }\hat e-\hat h=\hat0).\label{c2}
\end{eqnarray*}
\end{proof}

\newtheorem{Example21c}[Theorem2.1]{\bf Example}
\begin{Example21c}\label{Example21c}
\textnormal{Let $R=F[x,y,z]$, $k=x^3y^2z$ and let $B=\left[\begin{array}{cc}x^2y^2&x+1\\x^2&0\end{array}\right]$, $B'=\left[\begin{array}{cc}x^2y^2&0\\0&0\end{array}\right]$ and $
B''=\left[\begin{array}{cc}1+xyz&0\\0&0\end{array}\right]$. Note that $\widehat B$, $\widehat B'$ and $\widehat B''$ are $\langle x^3y^2z\rangle$-matrices.
Since $\textnormal{gcd}(x^2y^2,x^2)=x^2$ and $(\widehat{x+1})_{\langle{x^2}\rangle}$ is invertible in $R/\langle x^2\rangle$, it follows from Lemma \ref{Corollary2.5b}(ii) and Theorem \ref{Corollary12c}(a) that
\[\textnormal{Cen}(\widehat B_{\langle k\rangle})=\Theta(\textnormal{Cen}(B))=\left.\left\{\left[\begin{array}{cc}\hat a&(\widehat{x+1})\hat b\\\widehat{x^2}\hat b&\hat a+\widehat{x^2y^2}\hat b\end{array}\right]\right|\hat a,\hat b\in F[x,y,z]/\langle x^3y^2z\rangle\right\}.\]
Furthermore, it follows from Theorem \ref{Corollary12c}(b) that
\[\textnormal{Cen}(\widehat B'_{\langle k\rangle})=\left[\begin{array}{cc}R/\langle x^3y^2z\rangle&\langle \widehat{xz}\rangle\\\langle \widehat{xz}\rangle&R/\langle x^3y^2z\rangle\end{array}\right]\]
 and, since $\theta_{\langle x^3y^2z\rangle}(1+xyz)$ is invertible in $R/\langle x^3y^2z\rangle$, from Theorem \ref{Corollary12c}(c) that\begin{eqnarray*}\textnormal{Cen}(\widehat B'')&=&\Theta(\textnormal{Cen}(B''_{\langle k\rangle}))\ =\ \left[\begin{array}{cc}\textnormal{ann}(\hat f)\cap\textnormal{ann}(\hat g)&\textnormal{ann}(\hat g)\cap\textnormal{ann}(\hat e-\hat h)\\\textnormal{ann}(\hat f)\cap\textnormal{ann}(\hat e-\hat h)&\textnormal{ann}(\hat f)\cap\textnormal{ann}(\hat g)\end{array}\right]\\&=&\left[\begin{array}{cc}R/\langle x^3y^2z\rangle&\hat 0\\\hat 0&R/\langle x^3y^2z\rangle\end{array}\right].\end{eqnarray*}}
\end{Example21c}\vskip0.5cm

The following result is well-known.
\newtheorem{Lemma8a}[Theorem2.1]{\bf Lemma}
\begin{Lemma8a}\label{Lemma8a}
Let $R$ be a PID.
Then an element $\hat b\in R/\langle k\rangle$ is invertible if and only if~\textnormal{gcd}$(b,k)=1$.
\end{Lemma8a}

Using Lemma \ref{Lemma8a} and the fact that every matrix in $M_2(R)$ is a $\langle k\rangle$-matrix if $R$ is a PID, we simplify Theorem \ref{Corollary12c}(a) for the case when $R$ is a PID.

\newtheorem{PID}[Theorem2.1]{\bf Corollary}
\begin{PID}\label{PID}
Let $R$ be a PID and let $B=\left[\begin{array}{cc}e&f\\g&h\end{array}\right]\in M_2(R)$. Then
\[\textnormal{Cen}(\widehat B)=\Theta(\textnormal{Cen}(B))\]
if and only if $B$ is a scalar matrix or \textnormal{gcd}$(e-h,f,g,k)=1$.
\end{PID}

Note that although Corollary \ref{Corollary12b} is not a characterization of the $\langle k\rangle$-matrices for which (\ref{centr}) is true, it is easier to verify if Corollary \ref{Corollary12b} applies to a specific matrix in $M_2(R)$ than to verify if Theorem \ref{Corollary12c}(a) applies to a specific matrix in~$M_2(R)$.

\newtheorem{Corollary12b}[Theorem2.1]{\bf Corollary}
\begin{Corollary12b}\label{Corollary12b}
Let $R$ be a UFD, $k\in R$ and $B=\left[\begin{array}{cc}e&f\\g&h\end{array}\right]\in M_2(R)$. If at least one of the three elements $\hat e-\hat h$, $\hat f$ and $\hat g$ is invertible in $R/\langle k\rangle$, then
\begin{equation*}
\textnormal{Cen}(\widehat B)=\Theta(\textnormal{Cen}(B)).
\end{equation*}
\end{Corollary12b}

\begin{proof}
It follows trivially that $\widehat B$ is a $\langle k\rangle$-matrix. Without loss of generality, let us suppose that $\hat f$ is invertible in $R/\langle k\rangle$. Then, by Lemma \ref{Lemma8a} $\textnormal{gcd}(f,k)=1$. Hence condition (i) in Theorem \ref{Corollary12c}(a) is satisfied. Now, suppose that $\textnormal{gcd}(e-h,g,k)=\delta$. If $\delta$ is a unit, then condition (ii) is also satisfied. Thus suppose that $\delta$ is not a unit. Then, since $\hat f_{\langle k\rangle}$ is invertible in $R/\langle k\rangle$ and $\delta|k$, it follows that there is a $t\in R$ such that $tf=1+sk=1+sv\delta$ for some $s,v\in R$ which implies that $\hat t_{\langle \delta\rangle}\hat f_{\langle\delta\rangle}=\hat 1_{\langle\delta\rangle}$. Therefore condition (ii) in Theorem \ref{Corollary12c}(a) is also satisfied. 
\end{proof}

\section{The number of matrices in the centralizer of a matrix in $M_2(R/\langle k\rangle)$, $R$ a UFD and~$R/\langle k\rangle$ finite}\label{finite}

\noindent In this section $k\in R$ will always be a nonzero nonunit such that $R/\langle k\rangle$ is finite and we will always denote the number of elements in a ring $S$ by~$|S|$.
The purpose of this section is to determine the number of matrices in $\textnormal{Cen}_{M_2(R/\langle k\rangle)}(B)$, where $R$ is a UFD, $R/\langle k\rangle$ is finite and $B\in M_2(R/\langle k\rangle)$.

To reach our goal, we first need some preliminary results.

\newtheorem{Def4.2}{\bf Definition}[section]
\begin{Def4.2}\label{Def4.2}
Let $k\in R$, let $B=\left[\begin{array}{cc}e&f\\g&h\end{array}\right]\in M_2(R)$ and let $d:=\textnormal{gcd}(e-h,f,g,k)$. We define the relation~$\sim$ on $\textnormal{Cen}_{M_2(R/\langle k\rangle)}(\widehat B_{\langle k\rangle})$ as follows:
for $\widehat A_{\langle k\rangle}, \widehat C_{\langle k\rangle}\in \textnormal{Cen}_{M_2(R/\langle k\rangle)}(\widehat B_{\langle k\rangle})$,
\[\widehat A_{\langle k\rangle}\sim\widehat C_{\langle k\rangle}\qquad\textrm{iff}\qquad\widehat A_{\langle k\rangle}-\widehat C_{\langle k\rangle}\in M_2(\langle \widehat{(kd^{-1})}_{\langle k\rangle}\rangle).\]
\end{Def4.2}
It follows immediately that $\sim$ is an equivalence relation.

\noindent We denote the equivalence class of $\widehat A_{\langle k\rangle}$ by $\widehat A_{\langle k\rangle}^*$ and the set
\[\{\widehat A_{\langle k\rangle}^*\ |\ \widehat A_{\langle k\rangle}\in (\textnormal{Cen}_{M_2(R/\langle k\rangle)}(\widehat B_{\langle k\rangle}))\}\] of all equivalence classes by
\[(\textnormal{Cen}_{M_2(R/\langle k\rangle)}(\widehat B_{\langle k\rangle}))^*.\]

\noindent Since \[M_2(\langle\widehat{(kd^{-1})}_{\langle k\rangle}\rangle)\subseteq\left[\begin{array}{cc}\textnormal{ann}(\hat f_{\langle k\rangle})\cap\textnormal{ann}(\hat g_{\langle k\rangle})&\textnormal{ann}(\hat e_{\langle k\rangle}-\hat h_{\langle k\rangle})\cap \textnormal{ann}(\hat g_{\langle k\rangle})\\\textnormal{ann}(\hat e_{\langle k\rangle}-\hat h_{\langle k\rangle})\cap \textnormal{ann}(\hat f_{\langle k\rangle})&\textnormal{ann}(\hat f_{\langle k\rangle})\cap\textnormal{ann}(\hat g_{\langle k\rangle})\end{array}\right],\] it follows from Theorem \ref{Theorem21} that $M_2(\langle\widehat{(kd^{-1})}_{\langle k\rangle}\rangle)\subseteq \textnormal{Cen}_{M_2(R/\langle k\rangle)}(\widehat B_{\langle k\rangle})$. Therefore each equivalence class in $(\textnormal{Cen}_{M_2(R/\langle k\rangle)}(\widehat B_{\langle k\rangle}))^*$ has $|\langle \widehat{(kd^{-1})}_{\langle k\rangle}\rangle|^4$ elements.

\noindent We define addition $\boxplus$ and multiplication $\boxdot$ on $(\textnormal{Cen}_{M_2(R/\langle k\rangle)}(\widehat B_{\langle k\rangle}))^*$ by 
\begin{equation}\label{add}\widehat A_{\langle k\rangle}^*\boxplus\widehat C_{\langle k\rangle}^*=(\widehat A_{\langle k\rangle}+\widehat C_{\langle k\rangle})^*\quad\textnormal{ and }\quad\widehat A_{\langle k\rangle}^*\boxdot\widehat C_{\langle k\rangle}^*=(\widehat A_{\langle k\rangle}\boxdot\widehat C_{\langle k\rangle})^*.\end{equation}
 
\noindent It is easy to show that the binary operations $\boxplus$ and $\boxdot$ are well-defined and that the triple $\langle (\textnormal{Cen}_{M_2(R/\langle k\rangle)}(\widehat B_{\langle k\rangle}))^*,\boxplus,\boxdot\rangle$ is a ring, which we sometimes, if the context is clear, denote by $(\textnormal{Cen}_{M_2(R/\langle k\rangle)}(\widehat B_{\langle k\rangle}))^*.$ 

Using the following well-known result, Corollary \ref{direktesom} can easily be proved.

\newtheorem{finite4}[Def4.2]{\bf Theorem}
\begin{finite4}\label{finite4}
If $A_1,\ldots,A_m$ are ideals in a ring $S$ (not necessarily commutative or with a unit), then there is a monomorphism of rings
$\phi:S/(A_1\cap\cdots\cap A_m)\to S/A_1\oplus\cdots\oplus S/A_m$
defined by $\phi(s+(A_1\cap\cdots\cap A_m))=(s+A_1,\ldots,s+A_m).$ If $S^2+A_i=S$ for all $i$ and $A_i+A_j=S$ for all $i\neq j$, then $\phi$ is an isomorphism of rings.
\end{finite4}

\newtheorem{direktesom}[Def4.2]{\bf Corollary}
\begin{direktesom}\label{direktesom}
Let $R/\langle k\rangle$ be finite, and let $k=p_1^{n_1}p_2^{n_2}\cdots p_m^{n_m}$, with $p_1,\ldots,p_m$ different primes and $n_1,\ldots, n_m\ge1$. Then
\begin{itemize}
\item[\textnormal{(i)}]$\phi:R/\langle k\rangle\to R/\langle p_1^{n_1}\rangle\oplus R/\langle p_2^{n_2}\rangle\oplus\cdots \oplus R/\langle p_m^{n_m}\rangle$
defined by
\[\phi(\hat r)=(\theta_{\langle p_1^{n_1}\rangle}(r),\theta_{\langle p_2^{n_2}\rangle}(r),\cdots,\theta_{\langle p_m^{n_m}\rangle}(r))\]
is an isomorphism.
\item[\textnormal{(ii)}]$\Phi:M_2(R/\langle k\rangle)\to M_2(R/\langle p_1^{n_1}\rangle)\oplus\cdots\oplus M_2(R/\langle p_m^{n_m}\rangle)$
defined by
\[\Phi([\hat b_{ij}])=(\Theta_{\langle p_1^{n_1}\rangle}([b_{ij}]),\ldots,\Theta_{\langle p_m^{n_m}\rangle}([b_{ij}]))\]
is an isomorphism.
\end{itemize}
\end{direktesom}

We need the following trivial results in the next Lemma \ref{Lemma4.3}.

\newtheorem{cen1}[Def4.2]{\bf Lemma}
\begin{cen1}\label{cen1}
Let $S,S_1,\ldots, S_m$ be rings, $s\in S$ and let $\Gamma:S\to S_1\oplus\cdots\oplus S_m$ defined by $\Gamma(s)=(s_1,\ldots,s_m)$ be an isomorphism. Then
$t\in \textnormal{Cen}_S(s)$ if and only if $t_i\in\textnormal{Cen}_{S_i}(s_i),$
for all $i$.
\end{cen1}

\newtheorem{finite1}[Def4.2]{\bf Lemma}
\begin{finite1}\label{finite1}
Let $R/\langle k\rangle$ be finite.
An element $\hat b\in R/\langle k\rangle$ is invertible if and only if~\textnormal{gcd}$(b,k)=\nolinebreak1$.
\end{finite1}

\newtheorem{Lemma4.3}[Def4.2]{\bf Lemma}
\begin{Lemma4.3}\label{Lemma4.3}
Let $B=\left[\begin{array}{cc} e&f\\g&h\end{array}\right]\in M_2(R)$ and let $k\in R$. If \textnormal{gcd}$(e-\nolinebreak h,f,g,k)\linebreak=1$, then
\[|\textnormal{Cen}(\widehat B_{\langle k\rangle})|=|R/\langle k\rangle|^2.\] 
\end{Lemma4.3}

\begin{proof}
Suppose $k=p_1^{n_1}p_2^{n_2}\cdots p_m^{n_m}$, where $p_1,\ldots,p_m$ are different primes and $n_i\ge1$ for all $i$. It follows from Lemma \ref{direktesom}(ii) and Lemma~\ref{cen1} that
\[\textnormal{Cen}_{M_2(R/\langle k\rangle)}(\widehat B_{\langle k\rangle})\cong \bigoplus_{i=1}^{m}\textnormal{Cen}_{M_2(R/\langle p_i^{n_i}\rangle)}(\widehat B_{\langle p_i^{n_i}\rangle}).\]
Therefore,
\[|\textnormal{Cen}_{M_2(R/\langle k\rangle)}(\widehat B_{\langle k\rangle})|=\prod_{i=1}^{m}|\textnormal{Cen}_{M_2(R/\langle p_i^{n_i}\rangle)}(\widehat B_{\langle p_i^{n_i}\rangle})|.\] If we can show that $|\textrm{Cen}_{M_2(R/\langle{p_i^{n_i}}\rangle)}(\widehat B_{\langle p_i^{n_i}\rangle})|=|R/\langle p_i^{n_i}\rangle|^2,$ for all $i$, it follows, again from Lemma \ref{direktesom}(ii) and Lemma \ref{cen1}, that
\[|\textrm{Cen}_{M_2(R/\langle k\rangle)}(\widehat B_{\langle k\rangle})|=\prod_{i=1}^{m}|R/\langle p_i^{n_i}\rangle|^2=|R/\langle k\rangle|^2.\]

\noindent Let $p_i$ be an arbitrary prime in the prime factorization of $k$. Since \textnormal{gcd}$(e-h,f,g,k)=1$, it follows that $p_i\nmid f$ or $p_i\nmid g$ or $p_i\nmid e-h$. Thus, by Lemma~\ref{finite1}, at least one of~$\hat f_{\langle p_i^{n_i}\rangle}$, $\hat g_{\langle p_i^{n_i}\rangle}$ or~$\hat e_{\langle p_i^{n_i}\rangle}-\hat h_{\langle p_i^{n_i}\rangle}$ is invertible in $R/\langle p_i^{n_i}\rangle$.

\noindent If $\hat f$ is invertible in $R/\langle p_i^{n_i}\rangle$ with inverse $\hat t$, say, then given that $\textnormal{gcd}(e-h,f,g,p_i^{n_i})\linebreak=1$, it follows from Corollary \ref{Corollary12b} and Lemma \ref{Corollary2.5b}(ii) that
\begin{eqnarray}
\textnormal{Cen}_{M_2(R/\langle p_i^{n_i}\rangle)}(\widehat B)&=&\textnormal{Cen}\left(\left[\begin{array}{cc}\hat e&\hat f\\\hat g&\hat h\end{array}\right]\right)=\textnormal{Cen}\left(\left[\begin{array}{cc}\hat t\hat e&\hat 1\\\hat t\hat g&\hat t\hat h\end{array}\right]\right)
\nonumber\\
&=&\left.\left\{\hat a\left[\begin{array}{cc}\hat 1&\hat 0\\\hat 0&\hat 1\end{array}\right]+\hat b\left[\begin{array}{cc}\hat 0&\hat 1\\\hat t\hat g&-\hat t(\hat e-\hat h)\end{array}\right]\right|\hat a, \hat b\in R/\langle p_i^{n_i}\rangle\right\}.\label{1aaaa}
\end{eqnarray}
It can be similarly shown that if $\hat g$ is invertible in $R/\langle p_i^{n_i}\rangle$ with inverse $\hat t$, say, then
\begin{equation}\textnormal{Cen}_{M_2(R/\langle p_i^{n_i}\rangle)}(\widehat B)=\left.\left\{\hat a\left[\begin{array}{cc}\hat 1&\hat 0\\\hat 0&\hat 1\end{array}\right]+\hat b\left[\begin{array}{cc}\hat 0&\hat t\hat f\\\hat 1&-\hat t(\hat e-\hat h)\end{array}\right]\right|\hat a,\hat b\in R/\langle p_i^{n_i}\rangle\right\};\label{2}\end{equation}
and if $\hat e-\hat h$ is invertible in $R/\langle p_i^{n_i}\rangle$ with inverse $\hat t$, say, then
\begin{equation}\textnormal{Cen}_{M_2(R/\langle p_i^{n_i}\rangle)}(\widehat B)=\left.\left\{\hat a\left[\begin{array}{cc}\hat 1&\hat 0\\\hat 0&\hat 1\end{array}\right]+\hat b\left[\begin{array}{cc}\hat 0&-\hat t\hat f\\-\hat t\hat g&\hat 1\end{array}\right]\right|\hat a,\hat b\in R/\langle p_i^{n_i}\rangle\right\}.\label{3}\end{equation}
It is easy to see that the number of elements in the sets in (\ref{1aaaa}), (\ref{2}) and (\ref{3}) are $|R/\langle p_i^{n_i}\rangle|^2$.

\end{proof}

\newtheorem{Lemma4.4}[Def4.2]{\bf Lemma}
\begin{Lemma4.4}\label{Lemma4.4}
Let $k\in R$, let $B=\left[\begin{array}{cc} e&f\\g&h\end{array}\right]\in M_2(R)$ and let \[B'=\left[\begin{array}{cc}{d^{-1}}(e-h)&d^{-1}f\\d^{-1}g&0\end{array}\right],\] where $d:=\textnormal{gcd}(e-h,f,g,k)$, then
\[(\textnormal{Cen}_{M_2( R/\langle k\rangle)}(\widehat B_{\langle k\rangle}))^*\cong \textnormal{Cen}_{M_2(R/\langle kd^{-1}\rangle)}(\widehat{B'}_{\langle kd^{-1}\rangle}).\]

\end{Lemma4.4}
\begin{proof}
Since
\begin{eqnarray*}
&&\widehat A_{\langle k\rangle}^*\in (\textnormal{Cen}_{M_2(R/\langle k\rangle)}(\widehat B_{\langle k\rangle}))^*
\Leftrightarrow \widehat A_{\langle k\rangle}\in \textnormal{Cen}_{M_2(R/\langle k\rangle)}(\widehat B_{\langle k\rangle})\\
&\Leftrightarrow& \widehat A_{\langle k\rangle}\in \textnormal{Cen}_{M_2(R/\langle k\rangle)}\left(\left[\begin{array}{cc}\hat e_{\langle k\rangle}-\hat h_{\langle k\rangle}&\hat f_{\langle k\rangle}\\\hat g_{\langle k\rangle}&\hat 0_{\langle k\rangle}\end{array}\right]\right)\\
&\Leftrightarrow& A\left[\begin{array}{cc} e- h& f\\ g& 0\end{array}\right]-\left[\begin{array}{cc} e- h& f\\ g& 0\end{array}\right]A\in M_2(\langle k\rangle)\\
&\Leftrightarrow& AB'-B'A\in M_2(\langle kd^{-1}\rangle)\Leftrightarrow\widehat A_{\langle kd^{-1}\rangle}\in \textnormal{Cen}_{M_2(R/\langle kd^{-1}\rangle)}(\widehat{B'}_{\langle kd^{-1}\rangle})
\end{eqnarray*}
and
\begin{eqnarray*}\widehat A_{\langle k\rangle}^*=\widehat C_{\langle k\rangle}^*
&\Leftrightarrow&\widehat A_{\langle k\rangle}-\widehat C_{\langle k\rangle}\in M_2(\langle \widehat{kd^{-1}}_{\langle k\rangle}\rangle)\\&\Leftrightarrow& A-C\in M_2(\langle kd^{-1}\rangle)\Leftrightarrow \widehat A_{\langle kd^{-1}\rangle}=\widehat C_{\langle kd^{-1}\rangle}.\end{eqnarray*}
it follows that $\Gamma:(\textnormal{Cen}_{M_2( R/\langle k\rangle)}(\widehat B_{\langle k\rangle}))^*\to \textnormal{Cen}_{M_2(R/\langle kd^{-1}\rangle)}(\widehat{B'}_{\langle kd^{-1}\rangle})$, defined by \[\Gamma(\widehat A^*)=\widehat A_{\langle kd^{-1}\rangle},\] is a well-defined function which is $1-1$ and onto. It can be easily shown that $\Gamma$ is a homomorphism.
\end{proof}

We are finally able to determine the number of elements in the centralizer of a matrix in $M_2(R/\langle k\rangle),$ if $R$ is a UFD and $R/\langle k\rangle$ is finite.

\newtheorem{Theorem4.5}[Def4.2]{\bf Theorem}
\begin{Theorem4.5}\label{Theorem4.5}
 Suppose $R$ is a UFD, $k\in R$ is a nonzero nonunit such that $R/\langle k\rangle$ is finite, and \[B=\left[\begin{array}{cc}e&f\\g&h\end{array}\right]\in M_2(R),\] then \[|\textnormal{Cen}_{M_2(R/\langle k\rangle)}(\widehat B_{\langle k\rangle})|=|R/\langle kd^{-1}\rangle|^2\cdot|\langle (\widehat{kd^{-1}})_{\langle k\rangle}\rangle|^4,\]
where $d:=$\textnormal{gcd}$(e-h,f,g,k)$.
\end{Theorem4.5}
\begin{proof}
With $B'$ as in Lemma \ref{Lemma4.4}, it follows from Lemma \ref{Lemma4.3} that
\[|\textnormal{Cen}_{M_2(R/\langle kd^{-1}\rangle)}(\widehat {B'}_{\langle kd^{-1}\rangle})|=|R/\langle kd^{-1}\rangle|^2.\]
Since each equivalence class in $(\textnormal{Cen}_{M_2(R/\langle{k}\rangle)}(\widehat B_{\langle k\rangle}))^*$ has cardinality $|\langle (\widehat{kd^{-1}})_{\langle k\rangle}\rangle|^4$, it follows that
\begin{eqnarray*}
|\textnormal{Cen}_{M_2(R/\langle{k}\rangle)}(\widehat B_{\langle k\rangle})|&=&|(\textnormal{Cen}_{M_2(R/\langle{k}\rangle)}(\widehat B_{\langle k\rangle}))^*||\langle(\widehat{kd^{-1}})_{\langle k\rangle}\rangle|^4,\end{eqnarray*}
and so Lemma \ref{Lemma4.4} implies that
\begin{eqnarray*}|\textnormal{Cen}_{M_2(R/\langle{k}\rangle)}(\widehat B_{\langle k\rangle})|&=&|\textnormal{Cen}_{M_2(R/\langle{kd^{-1}}\rangle)}(\widehat{B'}_{\langle kd^{-1}\rangle})||\langle (\widehat{kd^{-1}})_{\langle k\rangle}\rangle|^4\\&=&|R/\langle kd^{-1}\rangle|^2|\langle (\widehat{kd^{-1}})_{\langle k\rangle}\rangle|^4.
\end{eqnarray*}
\end{proof}

\newtheorem{Example4.6}[Def4.2]{\bf Example}
\begin{Example4.6}\label{Example4.6}
\textnormal{
Let $R=\mathbb Z[i]$, $k=12$ so that $R/\langle k\rangle=\nolinebreak\mathbb Z_{12}[i]$ (see \cite{gauss}, p.\,604, Theorem 1) and let \[\widehat B=\left[\begin{array}{cc}\widehat{4i}&\hat 3+\widehat{6i}\\\widehat{9i}&\widehat{i}\end{array}\right].\]
Using the fact that every matrix is a $\langle k\rangle$-matrix if $R$ is a PID, note that, according to Lemma \ref{Corollary2.5b}(ii) and Theorem \ref{Theorem21}
\begin{eqnarray*}
\textnormal{Cen}_{M_2(\mathbb Z_{12}[i])}(\widehat B_{\langle 12\rangle})&=&\Theta_{\langle 12\rangle}\left(\left\{\left.\left[\begin{array}{cc}a&(1+2i)b\\3ib&a-3ib\end{array}\right]\right|a,b\in\mathbb Z[i]\right\}\right)+\left[\begin{array}{cc}\langle\hat 4\rangle&\langle \hat 4\rangle\\\langle\hat 4\rangle&\hat 0\end{array}\right]\\
&=&\left.\left\{\left[\begin{array}{cc}\hat a+\widehat{4c}&(\hat 1+\widehat{2i})\hat b+\widehat{4m}\\\widehat{3ib}+\widehat{4n}&\hat a-\widehat{3ib}\end{array}\right]\right|\hat a,\hat b, \hat c, \hat m, \hat n\in\mathbb{Z}_{12}[i]\right\}.
\end{eqnarray*}
Now, since $\textnormal{gcd}(3i, 3+6i,9i,12)=3$, let $d=3$ so that $kd^{-1}=12\cdot3^{-1}=4$. Since\[|\mathbb Z[i]/\langle 4\rangle|=|\{a+ib\ |\ a,b\in\mathbb Z_4\}|=16\qquad\textrm{and}\qquad|\langle\hat 4_{\langle 12\rangle}\rangle|=9\] it follows from Theorem \ref{Theorem4.5} that
\[|\textnormal{Cen}_{M_2(\mathbb Z_{12}[i])}(\widehat B_{\langle 12\rangle})|=16^2\cdot9^4=1679616.\]}
\end{Example4.6}

For $2\times 2$ matrices over a factor ring of $\mathbb Z$ we have the following result.
\newtheorem{Corollary4.5}[Def4.2]{\bf Corollary}
\begin{Corollary4.5}\label{Corollary4.5}
Let $\widehat B=\left[\begin{array}{cc} \hat e& \hat f\\ \hat g& \hat h\end{array}\right]\in M_2(\mathbb{Z}_k)$, then $|\textnormal{Cen}(\widehat B)|=(kd)^2,$
where $d=\textnormal{gcd}(e-h,f,g,k)$.
\end{Corollary4.5}
\begin{proof}
According to Theorem \ref{Theorem4.5}
\begin{equation*}
|\textnormal{Cen}_{M_2(\mathbb Z_k)}(\widehat B_{\langle k\rangle})|=|\mathbb Z_{kd^{-1}}|^2|\langle (\widehat{kd^{-1}})_{\langle k\rangle}\rangle|^4=(kd^{-1})^2d^4=(kd)^2.
\end{equation*}
\end{proof}
\newtheorem{Example4.5}[Def4.2]{\bf Example}
\begin{Example4.5}\label{Example4.5}
\textnormal{Let $\widehat B_{\langle 12\rangle}=\left[\begin{array}{cc}\hat 2_{\langle 12\rangle}&\hat 2_{\langle 12\rangle}\\\hat 4_{\langle 12\rangle}&\hat 8_{\langle 12\rangle}\end{array}\right].$ Since \textnormal{gcd}$(6,2,4,12)=2$, it follows that
\[|\textnormal{Cen}_{M_2(\mathbb Z_{12})}(\widehat B_{\langle 12\rangle})|=(12\cdot2)^2=24^2=576.\]}
\end{Example4.5}
\newtheorem{Question}[Def4.2]{\bf Remark}
\begin{Question}\label{Q1}
\textnormal{A natural example to include in this section, if such an example exists, would be one of a UFD $R$, which is not a PID, and a nonzero nonunit $k\in R$, such that $R/\langle k\rangle$ is finite. Unfortunately we could not find such an example. Neither have we been able to prove that if $R$ is UFD and $k\in R$ is a nonzero nonunit such that $R/\langle k\rangle$ is finite, then $R$ is a PID.}
\end{Question}

\bigskip

\end{document}